\documentclass[14pt]{article}

\usepackage[left=2.5cm, top=2cm, right=1cm, bottom=20mm, nohead, nofoot]{geometry}
\usepackage{cmap}
\usepackage[T2A]{fontenc}
\usepackage[utf8]{inputenc}
\usepackage[english]{babel}
\usepackage{amssymb}
\usepackage{graphicx}
\graphicspath{ {./img/} }
\usepackage{stmaryrd}
\usepackage{dsfont}
\usepackage{alltt}
\usepackage{caption}
\usepackage{comment}
\usepackage{amsmath,amsthm}
\usepackage[colorlinks=true, linkcolor=black, citecolor=black, urlcolor=black]{hyperref}
\usepackage{xcolor}

\newtheorem{theorem}{Theorem}[section]
\newtheorem{pro}[theorem]{Proposition}

\newtheorem{lemma}[theorem]{Lemma}
\theoremstyle{definition}

\newtheorem{remark}[theorem]{Remark}

\title{Extinction conditions for a two-sex branching process with random
mating and overlapping generations \\
N. Gordeev, M. Mota, E. Yarovaya

\thanks{\textbf{Nikita V. Gordeev}, Department of Probability Theory, Faculty of Mechanics and Mathematics, Lomonosov Moscow State University, Moscow, 119234, Russia, e-mail: nikita.gordeev@math.msu.ru.

\textbf{Corresponding author: Manuel Mota}, Department of Mathematics, University of Extremadura; and Institute of Advanced Scientific
Computation, University of Extremadura, Badajoz, Spain, e-mail: mota@unex.es. https://orcid.org/0000-0001-7012-9027

\textbf{Elena B. Yarovaya}, Department of Probability Theory, Faculty of Mechanics and Mathematics, Lomonosov Moscow State University, Moscow, 119234, Russia, e-mail: elena.yarovaya@math.msu.ru. https://orcid.org/0000-0002-6615-4315

}}

\date{}

\begin{document}
{\Large
\maketitle
\newpage
\tableofcontents
\newpage
\section{Introduction}

Two-sex branching processes represent an important class of mathematical
models describing the evolution of biological populations in which
individuals are divided into males and females. Unlike classical branching
processes, where each individual reproduces independently, in two-sex models
reproduction is possible only through the formation of mating pairs. This
significantly complicates the analysis, since it is necessary to take into
account not only the randomness in the number of offspring, but also the
randomness in pair formation, as well as the dependence of these processes
on the current sex-age structure of the population.

The theoretical foundations of two-sex branching processes were laid in the
work of Daley \cite{Daley1968}, who first introduced the two-sex
Galton-Watson process and obtained extinction conditions for the population.
Subsequently, the theory developed in various directions. In particular,
Molina, Mota and Ramos \cite{Molina2003} studied the two-sex Galton-Watson
process in varying environments, which made it possible to take into account
the influence of external factors on population dynamics. A comprehensive
survey of the literature on bisexual Galton-Watson branching processes can
be found in \cite{Hull2003}. In recent years, a number of generalizations
have been proposed that take into account random mating, dependence on
population size, and other biologically significant factors. Of particular
note is the work of Jacob, Molina and Mota \cite{Jacob2017}, which
introduced a general class of population-dependent two-sex processes with
random mating, and the recent paper of Molina and Mota \cite{Molina2024},
which for the first time introduced overlapping generations --- an important
characteristic of many biological species that reproduce throughout their
lives. A related recent contribution is due to Mart\'in-Gonz\'alez and
Gonz\'alez-Reyes \cite{MartinGonzalez2026}, who study the time to extinction
of a class of populations with sexual reproduction using mating functions
and connect branching processes with extreme value theory.

However, in the Molina and Mota \cite{Molina2024} model, the mortality
process is described quite simply: only individuals that have reached the
maximum age are removed from the population. In real populations, mortality
occurs at all age stages and may depend on population density, resource
availability, and other factors. This requires further generalization of
the model.

The present paper is devoted to the further development of the two-sex
branching process with random mating and overlapping generations. Unlike
the model of Molina and Mota \cite{Molina2024}, in which mortality is
trivial and all individuals survive up to the maximum age, we introduce a
more realistic mechanism of age-dependent survival, where the probabilities
of survival depend on the sex, age, and current state of the population.
To the best of our knowledge, such a generalization has not been previously
considered in the literature on two-sex branching processes with overlapping
generations. This assumption allows us to model populations in which
mortality affects all age groups and may be influenced by external factors
such as resource competition, predation, and population density.

We investigate the general properties of the proposed model and establish
sufficient conditions for both extinction and non-extinction of the
population. The obtained results generalize the theorem of Molina and
Mota \cite{Molina2024} to the case of nontrivial age-dependent survival
and provide a more biologically adequate tool for studying the long-term
dynamics of populations with overlapping generations.

The structure of the paper is as follows. Section 2 describes the class of
two-sex branching processes and the basic concepts related to the mating
function. Section 3 presents the definition of the two-sex branching
process with random mating and overlapping generations proposed by Molina
and Mota \cite{Molina2024} and introduces our modification with age-dependent
survival. Section 4 presents the general properties of the process,
including the absorbing nature of the zero state, transience of nonzero
states, and dichotomy. Section 5 proves the main limit theorems: a
sufficient condition for population extinction (Theorem 5.1) and a
sufficient condition for non-extinction (Theorem 5.2). Section 6 presents
the results of numerical simulations. The conclusion summarizes the results.

\newpage

\section{Class of two-sex branching processes}

A two-sex branching process models a population divided into males and
females, where reproduction is possible only through the formation of pairs.

Let $F_n$ and $M_n$ be the number of females and males at time $n$,
respectively. The number of mating pairs formed at time $n$ is determined
by the mating function
\[
L_n = L(F_n, M_n).
\]
For example, $L(F_n, M_n) = \min(F_n, M_n)$.

Each pair produces a random number of offspring
\[
\xi = (\xi^{(f)}, \xi^{(m)}),
\]
where $\xi^{(f)}$ and $\xi^{(m)}$ are respectively the number of females
and males among the offspring of one pair. Then
\[
F_{n+1} = \sum_{i=1}^{L_n} \xi_{i}^{(f)},
\qquad
M_{n+1} = \sum_{i=1}^{L_n} \xi_{i}^{(m)}.
\]

The main features of this class of processes are as follows:
\begin{itemize}
    \item the generation changes in two stages: pair formation and reproduction;
    \item the unit of reproduction is the pair, not an individual;
    \item the process is nonlinear, since $L_n$ depends on the current state
          $(F_n, M_n)$;
    \item the population may become extinct due to sex imbalance.
\end{itemize}

The scheme of random mating considered in this paper follows the approach
introduced by Jacob, Molina and Mota \cite{Jacob2017}, where the formation
of pairs is modelled through independent Bernoulli random variables rather
than through a deterministic mating function. This scheme is presented in
detail in Section 3.

\begin{remark}
Unlike the two-sex process, in the two-type process each individual
produces offspring independently of others. The key characteristic is the
possibility of type change at birth. A particle of type $i$ can give rise
to both particles of its own type and particles of the opposite type. The
change of generations occurs in one stage. Therefore, the class we consider
better describes complex organisms.
\end{remark}

\newpage

\section{Two-sex branching process with random mating and overlapping generations}

We begin by recalling the model of Molina and Mota \cite{Molina2024}, 
which describes a two-sex population with random mating and overlapping 
generations, in which individuals are removed only upon reaching the 
maximum age. For this model, we introduce the state space, the 
recurrence relations governing the evolution of the process, and the 
assumptions imposed on the random variables that describe mating, 
reproduction, and survival. We then propose a modification of this 
model, in which mortality is no longer trivial: the probability of 
surviving to the next time step depends on the sex, age, and current 
state of the population. For the modified process, we derive the 
corresponding recurrence relations and formulate the additional 
assumptions on the survival random variables.

\subsection{Basic model: random mating and overlapping generations}

Consider a biological population with sexual reproduction, in which females
and males reach a maximum age denoted by $a \ge 1$. Upon reaching this age,
they are removed from the population. During their lives, individuals can
reproduce, but before that they must form male-female pairs. The mating
probability of a female or male does not depend on age, but is determined
only by the current numbers of females and males in the population. However,
reproduction depends on the age of individuals, especially females.

Let $\mathbb{N}$ be the set of nonnegative integers and
\[
\mathbb{N}^{2a} :=
\underbrace{\mathbb{N} \times \cdots \times \mathbb{N}}_{2a}.
\]

The two-sex branching process with random mating and overlapping
generations, denoted by $\{X_n\}_{n \in \mathbb{N}}$, is defined as a
sequence of $2a$-dimensional random vectors
\[
X_n := \bigl(F_n^{(1)}, \ldots, F_n^{(a)}, \, M_n^{(1)}, \ldots, M_n^{(a)}\bigr),
\]
where $F_n^{(h)}$ and $M_n^{(h)}$ denote the number of females and males
respectively of age $h \in \{1,\ldots,a\}$ at time $n$.

The process is defined recursively. Initially,
\[
X_0 = \bigl(F_0^{(1)}, \ldots, F_0^{(a)}, \, M_0^{(1)}, \ldots, M_0^{(a)}\bigr) \in \mathbb{N}^{2a}.
\]
Given that at time $n$
\[
X_n = \bigl(F_n^{(1)}, \ldots, F_n^{(a)}, \, M_n^{(1)}, \ldots, M_n^{(a)}\bigr) \in \mathbb{N}^{2a},
\]
for each $h \in \{1,\ldots,a\}$ the offspring produced by females of age $h$
is given by
\[
\sum_{i=1}^{F_n^{(h)}} \delta_{n,i}^{(h)}(X_n)\, \xi_{n+1,i}^{(h)}(X_n),
\]
where:
\begin{itemize}
    \item $\delta^{(h)}_{n,i}(X_n)$ is a Bernoulli random variable taking
    value $1$ with probability $p(X_n)$ if the $i$-th female of age $h$
    forms a pair with a male at time $n$. Terms corresponding to pairs
    that were not formed are equal to $0$.

    \item $\xi^{(h)}_{n+1,i}(X_n)$ is a nonnegative two-dimensional random
    vector representing the number of females and males born at time
    $n+1$ by the pair in which the $i$-th female of age $h$ participates,
    provided that such a pair was indeed formed, i.e.
    $\delta^{(h)}_{n,i}(X_n) = 1$.
\end{itemize}

The total number of offspring born at time $n$ (having age $1$ at time
$n+1$) satisfies
\begin{equation}
(F^{(1)}_{n+1},\, M^{(1)}_{n+1}) =
\sum_{h=1}^{a}
\sum_{i=1}^{F^{(h)}_n}
\delta^{(h)}_{n,i}(X_n)\,
\xi^{(h)}_{n+1,i}(X_n).
\tag{3.1}
\end{equation}

The number of females and males of age $h+1$ at time $n+1$ equals the
number of individuals of age $h$ at time $n$, i.e.
\begin{equation}
(F^{(h+1)}_{n+1},\, M^{(h+1)}_{n+1}) = (F^{(h)}_n,\, M^{(h)}_n),
\qquad h = 1,\ldots,a-1.
\tag{3.2}
\end{equation}

Individuals of age $a$ at time $n$ are removed from the population.

Finally,
\[
X_{n+1} = (F^{(1)}_{n+1}, \ldots, F^{(a)}_{n+1},\,
M^{(1)}_{n+1}, \ldots, M^{(a)}_{n+1}).
\]

Obviously, if 
\[
\min\Bigl\{\sum_{h=1}^a F^{(h)}_n,\ \sum_{h=1}^a M^{(h)}_n\Bigr\} = 0,
\]
then
\[
\delta^{(h)}_{n,i}(X_n) = 0 \quad \text{almost surely}.
\]

We introduce the random vectors $\delta_n(x)$ and $\xi_{n+1}(x)$ governing
the mating and reproduction phases:
\[
\delta_n(x) :=
\bigl(\delta^{(h)}_{n,i}(x),\
h = 1,\ldots,a,\
i = 1,\ldots,F^{(h)}\bigr),
\]
\[
\xi_{n+1}(x) :=
\bigl(\xi^{(h)}_{n+1,i}(x),\
h = 1,\ldots,a,\
i = 1,\ldots,F^{(h)}\bigr).
\]

We assume the following conditions:

\begin{itemize}
    \item[(c1)] A female cannot mate with several males simultaneously, 
    unlike a male, who can mate with several females. Therefore, if we 
    denote by $\delta^{(h)}_{n,i}(x)$ the random variable representing 
    the total number of pairs formed by the $i$-th female of age $h$ at 
    time $n$, it is assumed that $\delta^{(h)}_{n,i}(x)$, for $h = 1,\ldots,a$ 
    and $i = 1,\ldots,F^{(h)}$, are independent and identically distributed 
    Bernoulli random variables with parameter $p(x)$, i.e.
    \[
    \mathbb{P}(\delta^{(h)}_{n,i}(x)=1) = p(x).
    \]

    \item[(c2)] The sequence $\{(\delta_n(x),\,\xi_{n+1}(x))\}_{n\in\mathbb{N}}$ 
    consists of independent and identically distributed random vectors.

    \item[(c3)] Given $\delta_n(x)$, the random vectors 
    $\xi^{(h)}_{n+1,i}(x)$ are independent and identically distributed, 
    and their distribution depends only on the value of $\delta^{(h)}_{n,i}(x)$. 
    More precisely,
    \[
    \mathbb{P}\bigl(\xi^{(h)}_{n+1,i}(x) = (0,0) \,\big|\, 
    \delta^{(h)}_{n,i}(x) = 0\bigr) = 1,
    \]
    and
    \[
    \mathbb{P}\bigl(\xi^{(h)}_{n+1,i}(x) = (k,l) \,\big|\, 
    \delta^{(h)}_{n,i}(x) = 1\bigr)
    = p^{(h)}_{k,l}(x), 
    \qquad k,l \in \mathbb{N},
    \]
    where $\{p^{(h)}_{k,l}(x)\}_{k,l\in\mathbb{N}}$ is a two-dimensional 
    probability distribution depending on the age $h$ of the female and 
    the population state $x$. 
    The conditional expectations  and variances are assumed to be finite for 
    any $x \in \mathbb{N}^{2a}$ and $h \in \{1,\ldots,a\}$ and are denoted by
    \[
    \varepsilon^{(h)}(x) = 
    \bigl(\varepsilon^{(h)}_f(x),\ \varepsilon^{(h)}_m(x)\bigr)
    :=
    \mathbb{E}\bigl[\xi^{(h)}_{n+1,i}(x) 
    \,\big|\, \delta^{(h)}_{n,i}(x)=1\bigr],
    \]
    \[
    \sigma^{(h)}_f(x):=
    {Var}\bigl[\xi^{(f,h)}_{n+1,i}(x) 
    \,\big|\, \delta^{(h)}_{n,i}(x)=1\bigr]\ ,\quad
    \sigma^{(h)}_m(x):=
    {Var}\bigl[\xi^{(m,h)}_{n+1,i}(x) 
    \,\big|\, \delta^{(h)}_{n,i}(x)=1\bigr],
    \]
    where $\varepsilon^{(h)}_f(x)$ and $\varepsilon^{(h)}_m(x)$ are the 
    conditional mean numbers of females and males born and $\sigma^{(h)}_f(x)$ and
    $\sigma^{(h)}_m(x)$ their corresponding variances.
\end{itemize}

It is clear that the two-sex process $\{X_n\}_{n\in\mathbb{N}}$ is a population-dependent multi-type branching process with $2a$ types, since there are several types of females and males in the population corresponding to all ages from $1$ to $a$.

\subsection{Modification: population-dependent survival condition}

In the model of Molina and Mota \cite{Molina2024}, only individuals that
have reached age $a$ die when transitioning to the next generation. We
add a survival condition that describes mortality within the population
depending on the sex and age of the individual and on the state of the
population at time $n$.

\begin{itemize}
    \item $\phi^{(h)}_{n,i}(X_n)$ is a Bernoulli random variable taking
    value $1$ with probability $q_f^{(h)}(X_n)$ if the $i$-th female of
    age $h$, living at time $n$, survives to time $n+1$.

    \item $\mu^{(h)}_{n,i}(X_n)$ is a Bernoulli random variable taking
    value $1$ with probability $q_m^{(h)}(X_n)$ if the $i$-th male of
    age $h$, living at time $n$, survives to time $n+1$.
\end{itemize}

Thus, conditions (3.1) and (3.2) take the form:
\begin{equation}
(F^{(1)}_{n+1},\, M^{(1)}_{n+1}) =
\sum_{h=1}^{a}
\sum_{i=1}^{F^{(h)}_n}
\delta^{(h)}_{n,i}(X_n)\,
\xi^{(h)}_{n+1,i}(X_n).
\tag{3.3}
\end{equation}
\begin{equation}
(F^{(h+1)}_{n+1},\, M^{(h+1)}_{n+1}) = 
\Bigl(\sum_{i=1}^{F^{(h)}_n} \phi^{(h)}_{n,i}(X_n),\
\sum_{i=1}^{M^{(h)}_n} \mu^{(h)}_{n,i}(X_n)\Bigr), 
\qquad h = 1,\ldots,a-1.
\tag{3.4}
\end{equation}

Obviously, $\phi^{(a)}_{n,i}(X_n) = \mu^{(a)}_{n,i}(X_n) = 0$ almost surely.

We introduce the random vectors $\phi_{n}(x)$ and $\mu_{n}(x)$ governing
the survival phase:
\[
\phi_n(x) := 
\bigl(\phi^{(h)}_{n,i}(x),\
h = 1,\ldots,a,\
i = 1,\ldots,F^{(h)}\bigr),
\]
\[
\mu_n(x) := 
\bigl(\mu^{(h)}_{n,i}(x),\
h = 1,\ldots,a,\
i = 1,\ldots,M^{(h)}\bigr).
\]

With respect to this sequence we formulate a new assumption:

\begin{itemize}
    \item[(c4)] The sequence $\{(\phi_n(x),\,\mu_{n}(x))\}_{n\in\mathbb{N}}$
    consists of independent and identically distributed random vectors and $\phi^{(h)}_{n,i}(x),\
h = 1,\ldots,a,\
i = 1,\ldots,F^{(h)}, \mu^{(h)}_{n,i}(x),\
h = 1,\ldots,a,\
i = 1,\ldots,M^{(h)}$ are independent random variables. Moreover the sequence $\{(\phi_n(x),\,\mu_{n}(x))\}_{n\in\mathbb{N}}$ is independent of the sequence $\{(\delta_n(x),\,\xi_{n+1}(x))\}_{n\in\mathbb{N}}$.
\end{itemize}

\section{General properties}

\begin{pro}
For the process $\{X_n\}_{n\in\mathbb{N}}$:

\begin{enumerate}
\item[(a)] The vector $\bar{0} \in \mathbb{N}^{2a}$ is an absorbing state.

\item[(b)] If for every $x \in \mathbb{N}^{2a}$, $x \neq \bar{0}$ at least one of the following conditions
holds:
    \begin{itemize}
        \item[(b1)] $p(x) < 1$,
        \item[(b2)] $p(x) > 0$ and 
        $\displaystyle \max\left\{ \sum_{k \in \mathbb{N}} p_{k,0}^{(h)}(x), 
        \sum_{l \in \mathbb{N}} p_{0,l}^{(h)}(x) \right\} > 0$ for
every $h \in \{1, \ldots, a\}$,
        \item[(b3)] $p(x)>0$, $q_f^{(1)}(x)<1$ and $\displaystyle  \sum_{l \in \mathbb{N}} p_{0,l}^{(1)}(x)  > 0$,
    \end{itemize}

then every nonzero state of the process is transient.
\end{enumerate}
\end{pro}

\begin{proof}
(a) It follows from the definition of the process $\{X_n\}_{n\in\mathbb{N}}$.

(b) The proof that (b1) or (b2) imply the transience of all the states of the process is the same as the proof of Proposition 3.2 in \cite{Molina2024}.

Let us prove that (b3) implies that every state $x=(f^{(1)}_x,\dots,f^{(a)}_x,m^{(1)}_x,\dots,m^{(a)}_x)\in \mathbb{N}^{2a}$ is transient. Denote $F_n=(F_n^{(1)},\dots,F_n^{(a)})$.

Since $\bar{0}$ is an absorbing state and $F_n=\bar{0}$ implies that $X_{n+1}=\bar{0}$, it will be enough to prove that $P(F_{n+a}=\bar{0}|X_n=x)>0$. 

First, we show that
\[
P(F_{n+a-1}^{(a)}=0,\dots,F_{n+a-1}^{(2)}=0\mid X_n=x)>0\tag{4.1}
\]
Denote $b_k(x)=P(F_{n+k}^{(k+1)}=0,\dots,F_{n+k}^{(2)}=0|X_n=x)$ and let us proceed by induction in $k$.

For $k=1$, considering the independence assumptions in (c4), we have that for every $x\in \mathbb{N}^{2a}$
$$
b_1(x):=P(F_{n+1}^{(2)}=0|X_n=x)=(1-q^{(1)}_f(x))^{f^{(1)}_x}>0
$$
Assume now that $b_1(x)>0, \dots, b_{k-1}(x)>0$. Denote, for $k=2,\dots,a$, $B_k=\{y\in \mathbb{N}^{2a}:f_y^{(2)}=0,\dots,f_y^{(k)}=0\}$. Using Markov property we have
\begin{align*}
b_k(x)\ge&\sum_{y\in B_k} P(F_{n+k}^{(k+1)}=0,\dots,F_{n+k}^{(2)}=0, X_{n+k-1}=y|X_n=x)\\
      =& \sum_{y\in B_k}P(F_{n+k}^{(2)}=0|X_{n+k-1}=y)P(X_{n+k-1}=y|X_n=x)\\
      =& \sum_{y\in B_k}b_1(y)P(X_{n+k-1}=y|X_n=x)
\end{align*}
By induction $b_1(y)>0$ for all $y\in B_k$. Also, by induction, $\sum_{y\in B_k}P(X_{n+k-1}=y|X_n=x)=b_{k-1}(x)>0$, so there exist $y\in B_k$ such that $P(X_{n+k-1}=y|X_n=x)>0$ and therefore $b_k(x)>0$. In particular, taking $k=a-1$, (4.1) is proved.

Now, we have that
\begin{align*}
P(F_{n+a}=\bar{0}|X_n=x)\ge&\sum_{y\in B_a}P(F_{n+a}=\bar{0},X_{n+a-1}=y|X_n=x)\\
=&\sum_{y\in B_a}P(F_{n+a}^{(2)}=0,F_{n+a}^{(1)}=0|X_{n+a-1}=y)P(X_{n+a-1}=y|X_n=x)\\
\end{align*}
By hypothesis, for every $y$, there exists $l$ such that $p_{0,l}^{(1)}(y)> 0$. Considering also the independence assumptions in (c3) and (c4),
\begin{align*}
&P(F_{n+a}^{(2)}=0,F_{n+a}^{(1)}=0|X_{n+a-1}=y)\\
=&P\left(\sum_{i=1}^{f_y^{(1)}}\phi_{n+a-1,i}^{(1)}(y)=0\right)P\left(\sum_{i=1}^{f_y^{(1)}}\delta_{n+a-1,i}^{(1)}(y)\xi_{n+a-1,i}^{(1,f)}(y)=0\right)\\
=&[(1-q_f^{(1)}(y))\left((1-p(y))+p(y)\sum_{l\in\mathbb{N}}p_{0,l}^{(1)}(y)\right)]^{f_y^{(1)}} > 0
\end{align*}
In addition, from (4.1), $\sum_{y\in B_a}P(X_{n+a-1}=y|X_n=x)=b_{a-1}(x)>0$, so there exists $y\in B_a$ such that $P(X_{n+a-1}=y|X_n=x)>0$ and therefore $P(F_{n+a}=\bar{0}|X_n=x)>0$ as we wanted to prove.
\end{proof}

\begin{remark}
Condition (b3) is specific of this model, because involves the probabilities of survival. Assuming that these probabilities are less than one for females of age 1, it is enough that the probability of not producing females is greater than 0  at this age to have a positive probability of leading the population to the extinction. The condition can be formulated also in terms of the males ($p(x)>0$, $q_m^{(1)}(x)<1$ and $ \sum_{l \in \mathbb{N}} p_{k,0}^{(1)}(x)  > 0$.
\end{remark}

\begin{theorem}[Molina and Mota \cite{Molina2024}]
Assume that for every nonzero state $x \in \mathbb{N}^{2a}$, $x \neq \bar{0}$,
condition (b) in Proposition 4.1 is satisfied. Then
\[
P\bigl( X_n \to \bar{0} \bigr) + P\bigl( \| X_n \| \to \infty \bigr) = 1,
\]
where $\| X_n \| = \sum_{h=1}^a ( F_n^{(h)} + M_n^{(h)} )$ is the total
population size at time $n$.
\end{theorem}

\begin{remark}
Since females play the main role in the evolution of the process (without females the population cannot reproduce, and unlimited growth of the number of females leads to unlimited growth of the entire population), it is of interest to establish the condition:
\[
P(F_n \to \bar{0}) + P(\|F_n\| \to \infty) = 1, \tag{4.2}
\]
where \( F_n := (F_n^{(1)}, \ldots, F_n^{(a)}) \) is the vector of female population sizes by age.
\end{remark}

\newpage

\section{Limit theorems}

\begin{theorem}[Sufficient condition for extinction]
Assume that for the process condition (4.2) holds and that
\[
\sum_{h=1}^a R_h \cdot \left( \prod_{j=1}^{h-1} S_j \right) < 1, \tag{5.1}
\]
where
\[
R_h := \sup_{x \in \mathbb{N}^{2a}} \bigl( p(x) \cdot \varepsilon_f^{(h)}(x) \bigr) > 0, \quad 
S_h := \sup_{x \in \mathbb{N}^{2a}} q_f^{(h)}(x), \quad h = 1,\ldots,a.
\]
Then the extinction probability of the population equals one:
\[
Q := P(X_n \to \bar{0}) = P(F_n \to \bar{0}) = 1.
\]
\end{theorem}

\begin{proof}
Consider the process with respect to females only. The vector 
$F_n = (F_n^{(1)}, \ldots, F_n^{(a)})$ evolves according to (3.3) and (3.4):
\[
F_{n+1}^{(1)} = \sum_{h=1}^a \sum_{i=1}^{F_n^{(h)}} \delta_{n,i}^{(h)}(X_n) \cdot \xi_{n+1,i}^{(f,h)}(X_n),
\]
\[
F_{n+1}^{(h+1)} = \sum_{i=1}^{F_n^{(h)}} \phi_{n,i}^{(h)}(X_n), \quad h = 1,\ldots,a-1,
\]
where $\xi_{n+1,i}^{(f,h)}(X_n)$ is the number of daughters born to the $i$-th female of age $h$.

For a fixed complete state 
\( x = (f^{(1)},\ldots,f^{(a)}, m^{(1)},\ldots,m^{(a)}) \in \mathbb{N}^{2a} \), 
the expectation of the vector \( F_{n+1} \) is:
\[
E[F_{n+1} \mid X_n = x] = f \cdot M_{11}(x),
\]
where \( f = (f^{(1)},\ldots,f^{(a)}) \) is the vector of female counts in state \( x \),
and \( M_{11}(x) \) is an \( a \times a \) matrix:
\[
M_{11}(x) = 
\begin{pmatrix}
p(x)\varepsilon_f^{(1)}(x) & q_f^{(1)}(x) & 0 & \cdots & 0 \\
p(x)\varepsilon_f^{(2)}(x) & 0 & q_f^{(2)}(x) & \cdots & 0 \\
\vdots & \vdots & \vdots & \ddots & \vdots \\
p(x)\varepsilon_f^{(a-1)}(x) & 0 & 0 & \cdots & q_f^{(a-1)}(x) \\
p(x)\varepsilon_f^{(a)}(x) & 0 & 0 & \cdots & 0
\end{pmatrix}.
\]

Define the majorizing matrix \(M_{11}^*\) of size \(a \times a\):
\[
M_{11}^* = 
\begin{pmatrix}
R_1 & S_1 & 0 & 0 & \cdots & 0 \\
R_2 & 0 & S_2 & 0 & \cdots & 0 \\
R_3 & 0 & 0 & S_3 & \cdots & 0 \\
\vdots & \vdots & \vdots & \vdots & \ddots & \vdots \\
R_{a-1} & 0 & 0 & 0 & \cdots & S_{a-1} \\
R_a & 0 & 0 & 0 & \cdots & 0
\end{pmatrix}.
\]

Then for any state \(x\) the componentwise inequality holds:
\[
(M_{11}(x))_{ij} \leq (M_{11}^*)_{ij}, \quad i,j = 1,\ldots,a.
\]

The spectral radius \(\rho\) of the matrix \(M_{11}^*\) is the unique positive root of the characteristic equation:
\[
\rho^a - \sum_{h=1}^a R_h \left( \prod_{j=1}^{h-1} S_j \right) \rho^{a-h} = 0. \tag{5.2}
\]

Consider the function \(f(\lambda) = \lambda^a - R_1 \lambda^{a-1} - R_2 S_1 \lambda^{a-2} - \cdots - R_a (S_1 S_2 \cdots S_{a-1})\). 
Note that \(f(0) = -R_a (S_1 S_2 \cdots S_{a-1}) < 0\) and \(\lim_{\lambda \to \infty} f(\lambda) = +\infty\). 
By continuity of \(f\) and uniqueness of the positive root, \(f(\lambda) < 0\) for \(0 \le \lambda < \rho\) and \(f(\lambda) > 0\) for \(\lambda > \rho\). 

Compute \(f(1) = 1 - \sum_{h=1}^a R_h \left( \prod_{j=1}^{h-1} S_j \right)\).

If condition (5.1) holds, i.e. \(\sum_{h=1}^a R_h \left( \prod_{j=1}^{h-1} S_j \right) < 1\), then \(f(1) > 0\). 
This means that \(1 > \rho\) (since at \(\lambda = 1\) the function is positive, which is possible only if \(1\) is to the right of the unique root). 
Therefore, \(\rho < 1\).

Let \(\nu = (\nu_1, \nu_2, \ldots, \nu_a)^\top\) be the positive right eigenvector of the matrix \(M_{11}^*\) corresponding to the eigenvalue \(\rho\):
\[
M_{11}^* \nu = \rho \nu.
\]

Define the process:
\[
\tilde{W}_n := \frac{F_n \nu}{\rho^n}.
\]

Then:
\[
\begin{aligned}
E[\tilde{W}_{n+1} \mid X_n = x]  
&= \frac{1}{\rho^{n+1}} E[F_{n+1} \nu \mid X_n = x]= \frac{1}{\rho^{n+1}} \bigl( F_n M_{11}(x) \bigr) \nu \leq \\
&\leq \frac{1}{\rho^{n+1}} F_n M_{11}^* \nu= \frac{1}{\rho^{n+1}} F_n (\rho \nu) = \frac{1}{\rho^n} F_n \nu = \tilde{W}_n.
\end{aligned}
\]

Thus, \(\{\tilde{W}_n\}\) is a nonnegative supermartingale with respect to the filtration \(\{\mathcal{F}_n\}\). By Doob's supermartingale convergence theorem, a nonnegative supermartingale converges almost surely to a finite limit. Therefore, there exists a random variable \(\tilde{W}_\infty\) such that:
\[
\tilde{W}_n \xrightarrow{n \to \infty} \tilde{W}_\infty \text{ a.s.}, \quad P(\tilde{W}_\infty < \infty) = 1.
\]

From the convergence of \(\tilde{W}_n\), it follows that:
\[
W_n = \rho^n \tilde{W}_n.
\]

Since \(\rho < 1\), we have \(\rho^n \to 0\). Consequently:
\[
W_n \xrightarrow{n \to \infty} 0 \text{ a.s.}.
\]

But \(W_n = \sum_{h=1}^a F_n^{(h)} \nu_h\), where all \(\nu_h > 0\). Therefore \(W_n \to 0\) almost surely implies \(F_n^{(h)} \to 0\) for all \(h\), i.e. \(F_n \to \bar{0}\) almost surely. Hence:
\[
P(\|F_n\| \to \infty) = 0.
\]

Then from condition (4.2) we obtain:
\[
P(F_n \to \bar{0}) = 1.
\]

If all females become extinct (\(F_n \to \bar{0}\)), then the entire population becomes extinct, since reproduction is impossible without females. Thus:
\[
Q = P(X_n \to \bar{0}) = 1.
\]
\end{proof}

\begin{remark}
The condition (5.1) \(\sum_{h=1}^a R_h \cdot \left( \prod_{j=1}^{h-1} S_j \right) < 1\) has a biological interpretation:
\begin{itemize}
    \item \(R_h\) is the maximum expected number of daughters from females of age \(h\);
    \item \(\prod_{j=1}^{h-1} S_j\) is the maximum probability of surviving to age \(h\);
    \item the product \(R_h \cdot \left( \prod_{j=1}^{h-1} S_j \right)\) is the maximum expected contribution of females to the population by age \(h\);
    \item the sum over all ages is the total maximum expected contribution of one female over her entire life to the replenishment of the population with females.
\end{itemize}
If this sum is less than one, then even under the most favorable conditions each female produces on average less than one daughter during her lifetime, which leads to inevitable extinction of the population.
\end{remark}

\newpage

In order to provide a sufficient condition for a positive probability of survival of the process, we will need the following lemma.

\begin{lemma}
\label{lemma}
Let $F_{n+1}=(F_{n+1}^{(1)},\dots,F_{n+1}^{(a)})$, 
$x=(f^{(1)},\dots,f^{(a)},m^{(1)},\dots,m^{(a)})\in \mathbb{N}^{2a}$ 
and $\nu=(\nu_1,\dots,\nu_a)\in \mathbb{N}^a$. Then
\begin{align*}
\operatorname{Var}[F_{n+1}\nu \mid X_n=x]
={}&\nu_1^2\sum_{h=1}^a f^{(h)}\bigl[(\varepsilon^{(h)}_f(x))^2 p(x)(1-p(x)) \\
&+ \sigma_f^{(h)}(x)\,p(x)\bigr]
+ \sum_{h=2}^{a}\nu_h^2 f^{(h-1)} q_f^{(h-1)}(x)\bigl(1-q_f^{(h-1)}(x)\bigr).
\end{align*}
\end{lemma}

\begin{proof}
The proof is a matter of straightforward calculation, taking into account 
the independence hypotheses assumed in (c3) and (c4).
\end{proof}

\begin{theorem}[Sufficient condition for non-extinction]
Assume that
\[
\sum_{h=1}^a r_h \cdot \left( \prod_{j=1}^{h-1} s_j \right) > 1, \tag{5.3}
\]
where
\[
r_h := \inf_{x \in \mathbb{N}^{2a}} \bigl( p(x) \cdot \varepsilon_f^{(h)}(x) \bigr) > 0, \qquad
s_h := \inf_{x \in \mathbb{N}^{2a}} q_f^{(h)}(x) > 0, \qquad h = 1, \ldots, a.
\]
Suppose in addition that, for every $h=1,\dots,a$, 
$\sup_x\varepsilon^{(h)}_f(x)$ and $\sup_x\sigma_f^{(h)}(x)$ are finite.
Then the population survives with positive probability:
\[
P(\|F_n\| \to \infty) > 0.
\]
\end{theorem}

\begin{proof}
Introduce the \emph{minoring} matrix $\check{M}_{11}$ of size $a \times a$, 
which is the natural counterpart of the majorizing matrix $M_{11}^*$:
\[
\check{M}_{11} :=
\begin{pmatrix}
r_1 & s_1 & 0 & \cdots & 0 \\
r_2 & 0 & s_2 & \cdots & 0 \\
\vdots & \vdots & \vdots & \ddots & \vdots \\
r_{a-1} & 0 & 0 & \cdots & s_{a-1} \\
r_a & 0 & 0 & \cdots & 0
\end{pmatrix}.
\]
By the definition of $r_h$ and $s_h$ as infima, we have the componentwise 
inequality
\[
\bigl(\check{M}_{11}\bigr)_{ij} \le \bigl(M_{11}(x)\bigr)_{ij}, 
\qquad i,j = 1, \ldots, a,
\]
for every state $x \in \mathbb{N}^{2a}$.

Let $\check{\rho}$ be the spectral radius of $\check{M}_{11}$. By the same 
argument as in the proof of Theorem 5.1, $\check{\rho}$ is the unique 
positive root of
\[
\check{\rho}^a - \sum_{h=1}^a r_h \left( \prod_{j=1}^{h-1} s_j \right) \check{\rho}^{a-h} = 0. \tag{5.5}
\]
Consider the function
\[
g(\lambda) := \lambda^a - \sum_{h=1}^a r_h \left( \prod_{j=1}^{h-1} s_j \right) \lambda^{a-h}.
\]
We have $g(0) = -r_a (s_1 \cdots s_{a-1}) < 0$ and $\lim_{\lambda \to \infty} g(\lambda) = +\infty$. 
Since $\check{\rho}$ is the unique positive root of $g$, we obtain
$g(\lambda) < 0$ for $0 \le \lambda < \check{\rho}$ and $g(\lambda) > 0$ 
for $\lambda > \check{\rho}$. By condition (5.3),
\[
g(1) = 1 - \sum_{h=1}^a r_h \left( \prod_{j=1}^{h-1} s_j \right) < 0.
\]
Since $g(\check{\rho}) = 0$ and $g$ is strictly increasing on 
$(\check{\rho}, +\infty)$, the inequality $g(1) < 0 = g(\check{\rho})$ 
implies
\[
1 < \check{\rho}.
\]

Let $\check{\nu} = (\check{\nu}_1, \ldots, \check{\nu}_a)^\top$ be the 
positive right eigenvector of $\check{M}_{11}$ corresponding to $\check{\rho}$:
\[
\check{M}_{11} \check{\nu} = \check{\rho} \check{\nu}, \qquad \check{\nu}_h > 0, 
\quad h = 1, \ldots, a.
\]
Such a vector exists by the Perron--Frobenius theorem, since $\check{M}_{11}$ 
is a nonnegative irreducible matrix. Set
\[
\check{\nu}_{\min} := \min_{1 \le h \le a} \check{\nu}_h > 0, \qquad
\check{\nu}_{\max} := \max_{1 \le h \le a} \check{\nu}_h.
\]
For every $n \ge 0$,
\[
\check{\nu}_{\min} \sum_{h=1}^a F_n^{(h)} 
\le F_n \check{\nu} 
= \sum_{h=1}^a F_n^{(h)} \check{\nu}_h 
\le \check{\nu}_{\max} \sum_{h=1}^a F_n^{(h)}.
\]
Consequently,
\[
\{ F_n \check{\nu} \to \infty \} = \bigl\{ \|F_n\| \to \infty \bigr\}. \tag{5.6}
\]

Choose $\varepsilon > 0$ such that
\[
\delta := \check{\rho} - \varepsilon > 1,
\]
and define for every $n \ge 0$ the event
\[
A_n := \bigl\{ F_{n+1} \check{\nu} > \delta \, F_n \check{\nu} \bigr\}.
\]
On $\bigcap_{n=0}^\infty A_n$ we have by iteration
\[
F_n \check{\nu} > \delta^n F_0 \check{\nu}.
\]
Since $\delta > 1$ and $F_0 \check{\nu} > 0$, we obtain 
$F_n \check{\nu} \to \infty$, and therefore by (5.6) 
$\|F_n\| \to \infty$. Thus
\[
\bigl\{ \|F_n\| \to \infty \bigr\} \supseteq \bigcap_{n=0}^\infty A_n,
\]
and consequently
\[
P(\|F_n\| \to \infty) \ge P\left( \bigcap_{n=0}^\infty A_n \right). \tag{5.7}
\]

It remains to show that $P(\bigcap_{n=0}^\infty A_n) > 0$. 

For a fixed state $x$ with $f = (f^{(1)}, \ldots, f^{(a)})$, we set 
$y_n := F_n \check{\nu}$ and $Y_{n+1} := F_{n+1} \check{\nu}$. Using 
the componentwise inequality $\check{M}_{11} \le M_{11}(x)$, we obtain
\[
E[Y_{n+1} \mid X_n = x] 
= \bigl( f \, M_{11}(x) \bigr) \check{\nu} 
\ge f \, \check{M}_{11} \check{\nu} 
= \check{\rho} \, f \check{\nu} 
= \check{\rho} \, y_n.
\]
Hence the conditional mean of $Y_{n+1}$ exceeds the threshold 
$(\check{\rho} - \varepsilon) y_n$ by at least $\varepsilon y_n$:
\[
E[Y_{n+1} \mid X_n = x] - (\check{\rho} - \varepsilon) y_n 
\ge \check{\rho} y_n - (\check{\rho} - \varepsilon) y_n 
= \varepsilon y_n > 0.
\]

The complement of $A_n$ has the form
\[
A_n^c = \bigl\{ Y_{n+1} \le (\check{\rho} - \varepsilon) y_n \bigr\}.
\]
For every outcome in $A_n^c$,
\[
Y_{n+1} - E[Y_{n+1} \mid X_n = x] 
\le (\check{\rho} - \varepsilon) y_n - E[Y_{n+1} \mid X_n = x] 
\le (\check{\rho} - \varepsilon) y_n - \check{\rho} y_n 
= -\varepsilon y_n.
\]
Therefore, as a set-theoretic inclusion,
\[
A_n^c \subseteq \bigl\{ Y_{n+1} - E[Y_{n+1} \mid X_n = x] \le -\varepsilon y_n \bigr\}
\subseteq \bigl\{ \bigl| Y_{n+1} - E[Y_{n+1} \mid X_n = x] \bigr| \ge \varepsilon y_n \bigr\}.
\]
Applying Chebyshev's inequality to the random variable $Y_{n+1}$ 
conditioned on $X_n = x$, with $t := \varepsilon y_n > 0$, we get
\[
P\bigl( \bigl| Y_{n+1} - E[Y_{n+1} \mid X_n = x] \bigr| 
\ge \varepsilon y_n \,\big|\, X_n = x \bigr) 
\le \frac{\operatorname{Var}[Y_{n+1} \mid X_n = x]}{(\varepsilon y_n)^2}.
\]
Combining the last two displays yields
\[
P(A_n^c \mid X_n = x) 
\le \frac{\operatorname{Var}[Y_{n+1} \mid X_n = x]}{(\varepsilon y_n)^2}. \tag{5.8}
\]

By Lemma~\ref{lemma}, taking into account that 
$\sup_x\varepsilon^{(h)}_f(x)$ and $\sup_x\sigma_f^{(h)}(x)$ are finite, 
there exists a constant $C_2 > 0$ such that
\[
\operatorname{Var}[Y_{n+1} \mid X_n = x] \le C_2 \sum_{h=1}^a f^{(h)} = C_2 \|f\|_1.
\]
Moreover, since $\check{\nu}_h \ge \check{\nu}_{\min} > 0$ for all $h$,
\[
y_n = F_n \check{\nu} \ge \check{\nu}_{\min} \|F_n\|_1.
\]
Substituting into (5.8), we obtain
\[
P(A_n^c \mid X_n = x) 
\le \frac{C_2 \|F_n\|_1}{\varepsilon^2 \check{\nu}_{\min}^2 \|F_n\|_1^2} 
= \frac{C_3}{\varepsilon^2 \|F_n\|_1}, \tag{5.9}
\]
where $C_3 := C_2 / \check{\nu}_{\min}^2$.

Set $B_{n-1} := \bigcap_{k=0}^{n-1} A_k$ for $n \ge 1$ and $B_{-1} := \Omega$. 
On $B_{n-1}$, iterating the definition of $A_k$ gives
\[
F_n \check{\nu} > \delta^n F_0 \check{\nu},
\]
hence
\[
\|F_n\|_1 \ge \frac{F_n \check{\nu}}{\check{\nu}_{\max}} 
> \frac{\delta^n F_0 \check{\nu}}{\check{\nu}_{\max}} 
=: c_0 \, \delta^n,
\]
where $c_0 := F_0 \check{\nu} / \check{\nu}_{\max} > 0$. Combining this bound 
with (5.9), we find
\[
P(A_n \mid B_{n-1}) = 1 - P(A_n^c \mid B_{n-1}) 
\ge 1 - \frac{C_3}{\varepsilon^2 c_0 \delta^n}.
\]

By the multiplication rule for conditional probabilities,
\[
P\left( \bigcap_{n=0}^N A_n \right) 
= \prod_{n=0}^N P(A_n \mid B_{n-1}) 
\ge \prod_{n=0}^N \left( 1 - \frac{C_3}{\varepsilon^2 c_0 \delta^n} \right).
\]
Letting $N \to \infty$,
\[
P\left( \bigcap_{n=0}^\infty A_n \right) 
\ge \prod_{n=0}^\infty \left( 1 - \frac{C_3}{\varepsilon^2 c_0 \delta^n} \right). \tag{5.10}
\]

Denote $a_n := C_3 / (\varepsilon^2 c_0 \delta^n)$. Since $\delta > 1$, 
the series
\[
\sum_{n=0}^\infty a_n 
= \frac{C_3}{\varepsilon^2 c_0} \cdot \frac{1}{1 - \delta^{-1}} 
< \infty
\]
converges. It is a standard fact that if $a_n \ge 0$ and 
$\sum_{n=0}^\infty a_n < \infty$, then the infinite product 
$\prod_{n=0}^\infty (1 - a_n)$ converges to a strictly positive limit, 
provided that $0 \le a_n < 1$ for all sufficiently large $n$. The latter 
holds because $a_n \to 0$ as $n \to \infty$. Hence the right-hand side 
of (5.10) is strictly positive, and therefore
\[
P\left( \bigcap_{n=0}^\infty A_n \right) > 0.
\]
Combining this with (5.7), we conclude
\[
P(\|F_n\| \to \infty) \ge P\left( \bigcap_{n=0}^\infty A_n \right) > 0.
\]
The theorem is proved.
\end{proof}

\begin{remark}
Condition (5.3) is the reverse inequality with respect to (5.1). Here
infima are used instead of suprema. The biological interpretation is as
follows: if even under the least favorable conditions each female produces
on average more than one daughter during her lifetime, then the population
survives with positive probability. Together with Theorem 5.1, this result
provides a complete description of the subcritical and supercritical
regimes of the process.
\end{remark}

\begin{remark}[Critical case and the trichotomy]
Theorems 5.1 and 5.2 provide sufficient conditions for extinction and 
non-extinction, respectively, but they do not cover the critical case, 
in which
\[
\lambda := \sum_{h=1}^a r_h \left( \prod_{j=1}^{h-1} s_j \right) 
\le 1 \le \sum_{h=1}^a R_h \left( \prod_{j=1}^{h-1} S_j \right) =: \Lambda.
\]
In the critical regime the long-term behaviour of the process depends on 
the precise form of the functions $p(x)$, $\varepsilon_f^{(h)}(x)$ and 
$q_f^{(h)}(x)$, and a finer analysis is required. 

In particular, one may expect that, in the critical case, the dichotomy 
established in condition (4.1) should be replaced by a trichotomy:
\[
P(F_n \to \bar{0}) + P(\|F_n\| \to \infty) + P(\|F_n\| \to c) = 1,
\]
where $c \in \mathbb{N}^{2a} \setminus \{\bar{0}\}$ is a non-degenerate 
finite limit. That is, in addition to extinction and unbounded growth, 
the population may stabilise at a finite positive level. A rigorous 
analysis of the critical case and of the corresponding trichotomy is a 
natural direction for further research and will be addressed in a 
forthcoming work.
\end{remark}

\newpage

\section{Numerical simulation}

To verify the theoretical results and visualize the behavior of the two-sex branching process, a computer simulation was developed in Python. The simulation allows one to observe the evolution of the age structure of the population over time and to analyze the influence of various parameters on population dynamics.

Let the maximum age of an individual be \(a = 110\) years. At each discrete time moment \(n\) (generation), the state of the population is described by the vectors

\[
F_n = (F_n^{(1)}, \ldots, F_n^{(a)}), \qquad M_n = (M_n^{(1)}, \ldots, M_n^{(a)}),
\]

where \(F_n^{(h)}\) and \(M_n^{(h)}\) are the numbers of females and males of age \(h \in \{1,\ldots,a\}\), respectively. Age \(h = 1\) corresponds to newborns, and individuals of age \(a\) are removed from the population after the completion of the evolutionary step.

At the initial time \(n = 0\), a uniform distribution of individuals by age is specified:

\[
F_0^{(h)} = C, \quad M_0^{(h)} = C, \qquad h = 1,\ldots,a,
\]

where \(C\) is the initial number of individuals of each age for both sexes.

Survival is modeled using the classical Gompertz-Makeham model. The mortality intensity at age \(x\) is given by

\[
\mu(x) = A + \alpha e^{\beta x},
\]

where \(A\) is the Makeham constant (background mortality independent of age), \(\alpha\) and \(\beta\) are the Gompertz parameters determining the scale and rate of mortality increase with age.

The survival function -- the probability of surviving to age \(x\) -- is obtained by integrating the mortality intensity:

\[
S(x) = \exp\left(-\int_0^x \mu(t)\,dt\right) = \exp\left(-Ax - \frac{\alpha}{\beta}(e^{\beta x} - 1)\right).
\]

For newborns \(S(0) = 1\), and as \(x \to \infty\) the function tends to zero. The probability for an individual of age \(h\) to survive to age \(h+1\) given that it survived to age \(h\) is

\[
p_{\text{surv}}(h) = \frac{S(h+1)}{S(h)}.
\]

This value is used in the simulation: for each female or male of age \(h\), a Bernoulli random variable is generated with success probability \(p_{\text{surv}}(h)\). If the result is one, the individual moves to the next age class; if zero, it dies.

Mating is modeled as follows. For each female of age \(h\), the probability of forming a pair with a male in the current generation is given by a constant \(p_{\text{mating}}\), the same for all ages. The number of pairs formed for each female is modeled as an independent Bernoulli random variable. If a female does not form a pair, she produces no offspring in that generation.

The distribution of the number of offspring is given by the multinomial law. Let \(M = 3\) be the maximum litter size. If a pair is formed, the number of born females \(k\) and males \(l\) follows the distribution

\[
p_{k,l} = \frac{M!}{k!\,l!\,(M - k - l)!}\, p_f^{\,k}\, p_m^{\,l}\, (1 - p_f - p_m)^{M - k - l},
\]

where \(p_f\) is the probability of giving birth to a female for one offspring, \(p_m\) is the probability of giving birth to a male, and \(k + l \le M\). The quantity \(1 - p_f - p_m\) corresponds to the probability that the offspring does not survive at birth or belongs to another type (within the simplification). All newborns enter the first age class (\(h = 1\)) of the next generation.

The evolution of the population from generation \(n\) to generation \(n+1\) occurs in two stages. In the first stage, for each age \(h\) and each female of this age, a pair is formed with probability \(p_{\text{mating}}\). Each formed pair produces a random number of females and males according to the multinomial distribution described above. All newborns are summed into \(F_{n+1}^{(1)}\) and \(M_{n+1}^{(1)}\). In the second stage, individuals age: for each age \(h = 1,\ldots,a-1\), each female of age \(h\) survives to age \(h+1\) with probability \(p_{\text{surv}}(h)\), and the same for males. Survivors join the corresponding age classes \(h+1\) in the next generation. Individuals that reach age \(a\) are removed from the population.

After the simulation is completed, three-dimensional scatter plots are constructed, where the horizontal axis represents the generation number \(n\), the vertical axis represents age \(h\), and the third axis represents the number of individuals. The plots display generations from \(n = 110\) to \(n = 150\), which allows observing the established regime of population dynamics. The color of the points corresponds to the type of plot: red -- females, blue -- males, green -- total population. The formatting of the abundance axis is done as integers with thousand separators for improved readability.

The simulation parameters are placed at the beginning of the program and can be easily changed. This allows investigating the influence of various factors on population dynamics: varying the mating probability, survival parameters, probabilities of birth of females and males, maximum litter size, and initial population size.
\newpage

\textbf{Result of numerical simulation:}
\begin{figure}[h]
\centering
\includegraphics[width=0.8\textwidth]{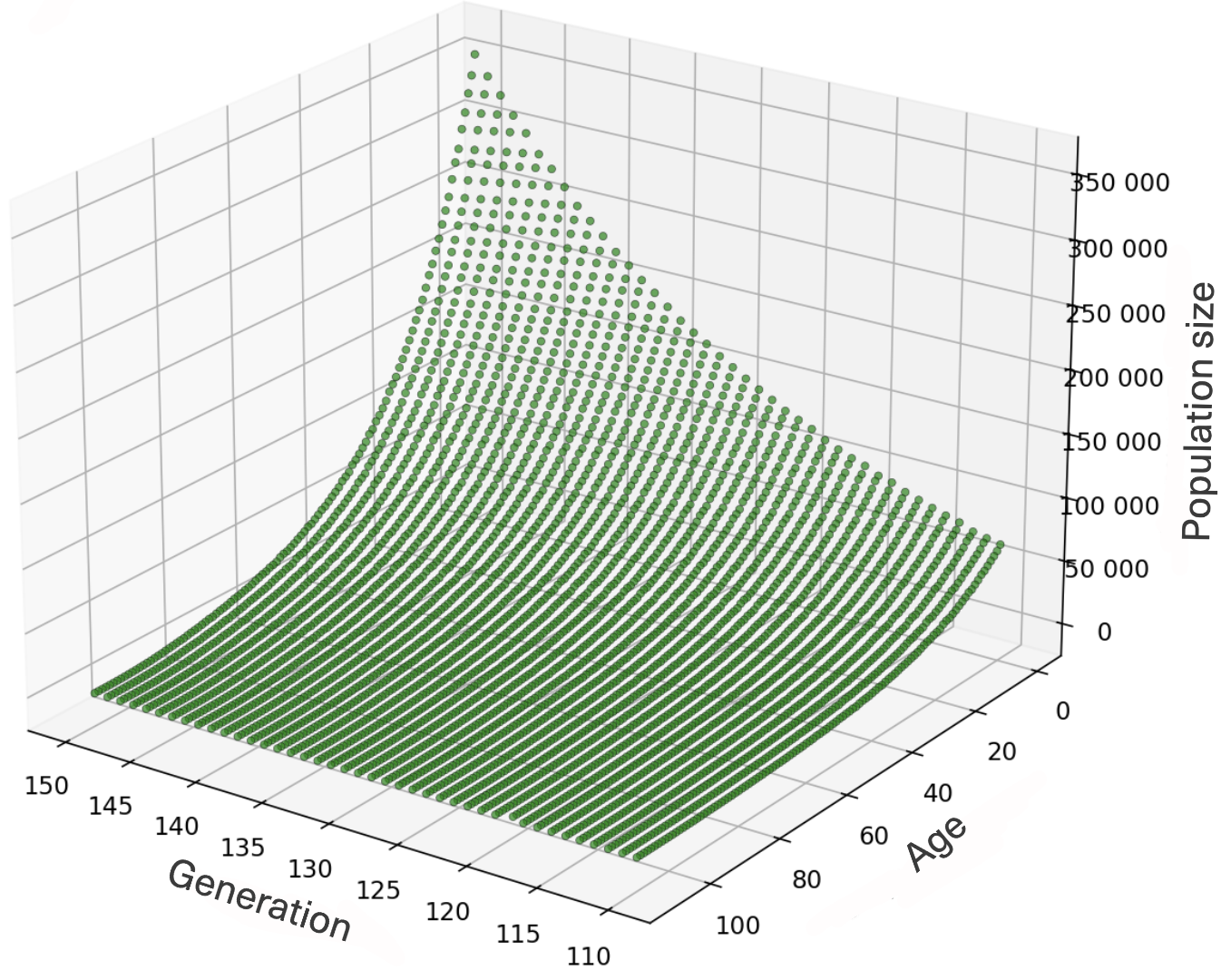}
\caption{Dynamics of the total population size by age and generation. 
Each point represents the total number of individuals (females and males 
combined) of a given age at a given generation. The horizontal axis shows 
the generation number $n$ (the discrete time index, i.e., the number of 
the current time step), the vertical axis shows age $h$, and the third 
axis shows the population size. The plot displays generations from $n = 110$ 
to $n = 150$ and ages from $h = 1$ to $h = 110$. The model parameters and 
initial conditions are listed below the figure.}
\label{fig:total}
\end{figure}

\vspace{0.3cm}
\noindent
\begin{minipage}[t]{0.48\textwidth}
\textit{Model parameters:}
\begin{itemize}
    \item $a = 110$ (max. age)
    \item $M = 3$ (max. litter size)
    \item $p_f = 0.3$, $p_m = 0.3$
    \item $A = 0.001$, $\alpha = 0.00001$, $\beta = 0.1$
    \item $p_{\text{mating}} = 0.05$
\end{itemize}
\end{minipage}
\hfill
\begin{minipage}[t]{0.48\textwidth}
\textit{Initial conditions:}
\begin{itemize}
    \item 100 individuals of each sex and age
\end{itemize}
\end{minipage}
\newpage
\noindent\textbf{Effect of mating probability on population dynamics:}

\begin{center}
\begin{minipage}{0.31\textwidth}
\centering
\includegraphics[width=\textwidth]{img/image1.png}
\vspace{0.2cm}
{ $p_{\text{mating}} = 0.05$}
\end{minipage}
\hfill
\begin{minipage}{0.31\textwidth}
\centering
\includegraphics[width=\textwidth]{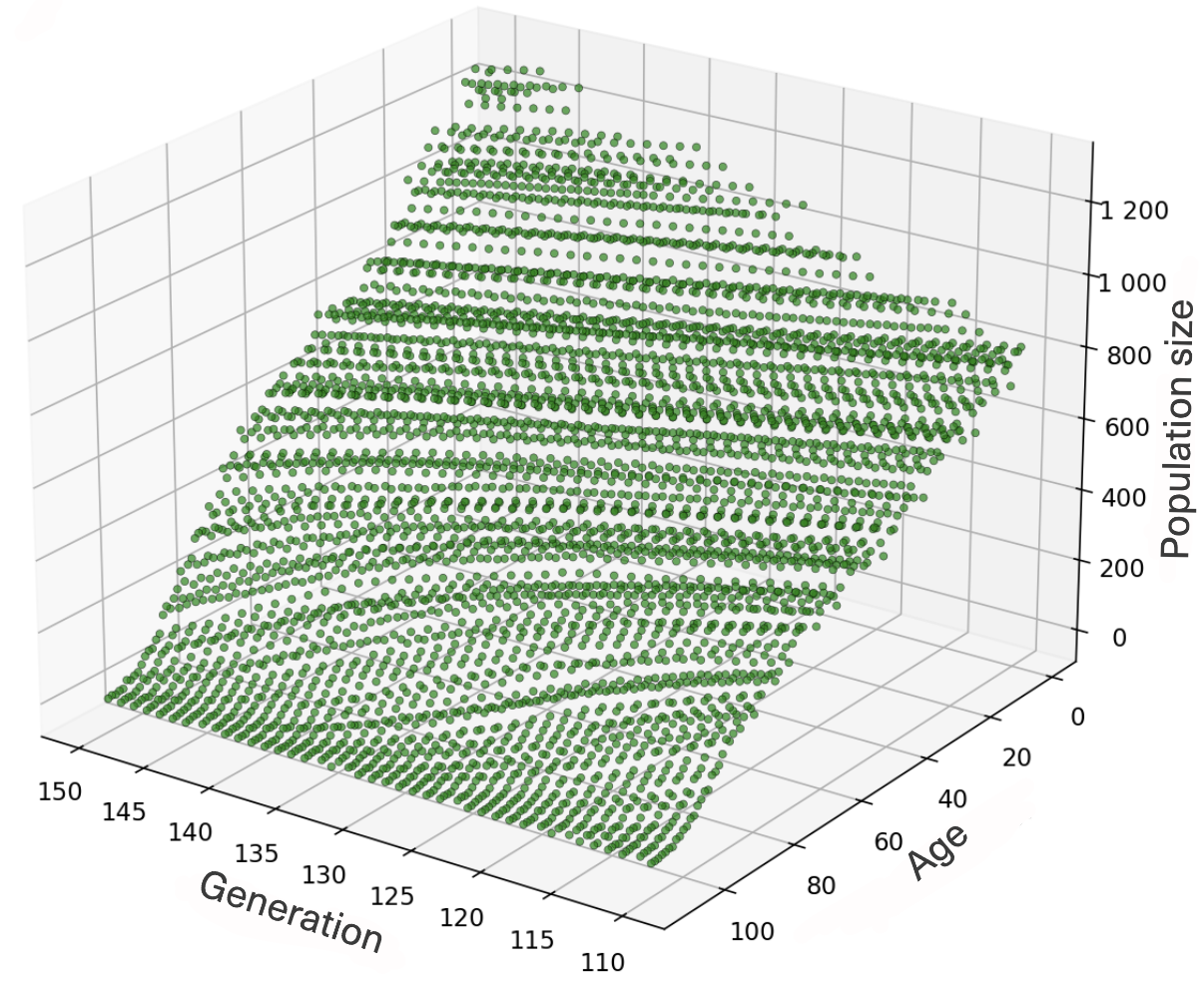}
\vspace{0.2cm}
{$p_{\text{mating}} = 0.02$}
\end{minipage}
\hfill
\begin{minipage}{0.31\textwidth}
\centering
\includegraphics[width=\textwidth]{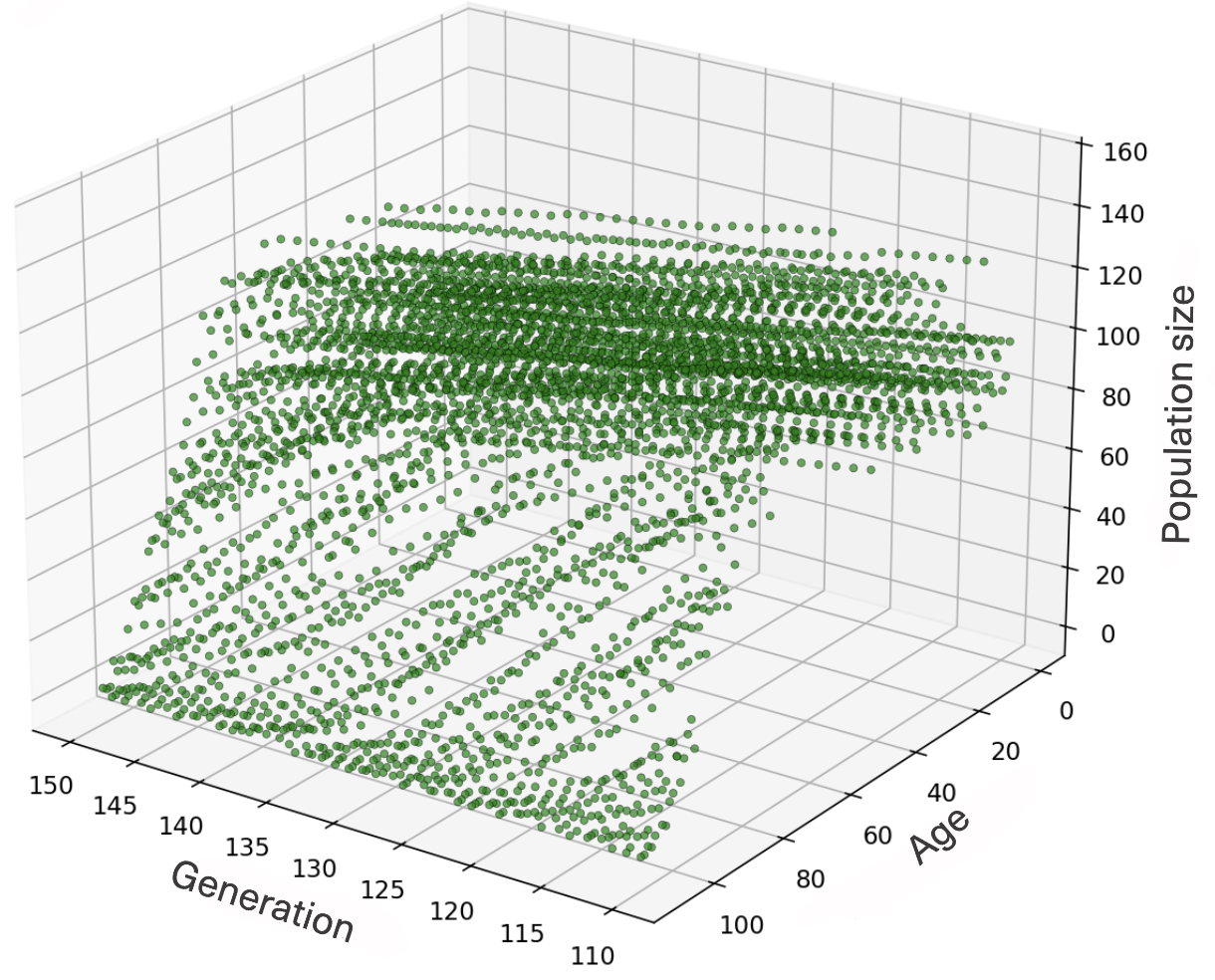}
\vspace{0.2cm}
{ $p_{\text{mating}} = 0.01$}
\end{minipage}
\end{center}

\noindent\textbf{Effect of mortality parameters on population dynamics:}

\begin{center}
\begin{minipage}{0.3\textwidth}
\centering
\includegraphics[width=\textwidth]{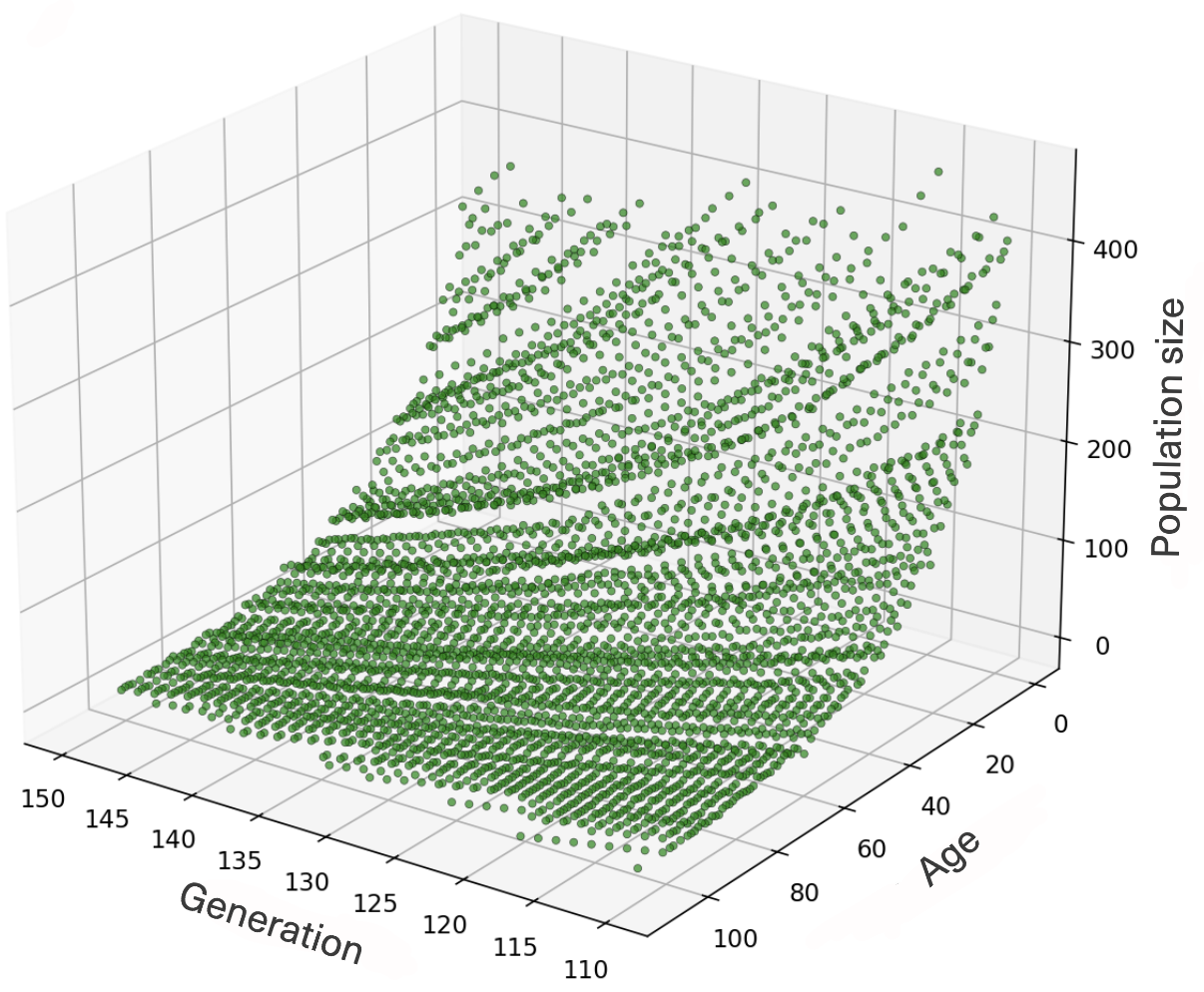}
\vspace{0.2cm}
$A = 0.05$ \\
$\alpha = 0.00001$ \\
$\beta = 0.1$
\end{minipage}
\hfill
\begin{minipage}{0.3\textwidth}
\centering
\includegraphics[width=\textwidth]{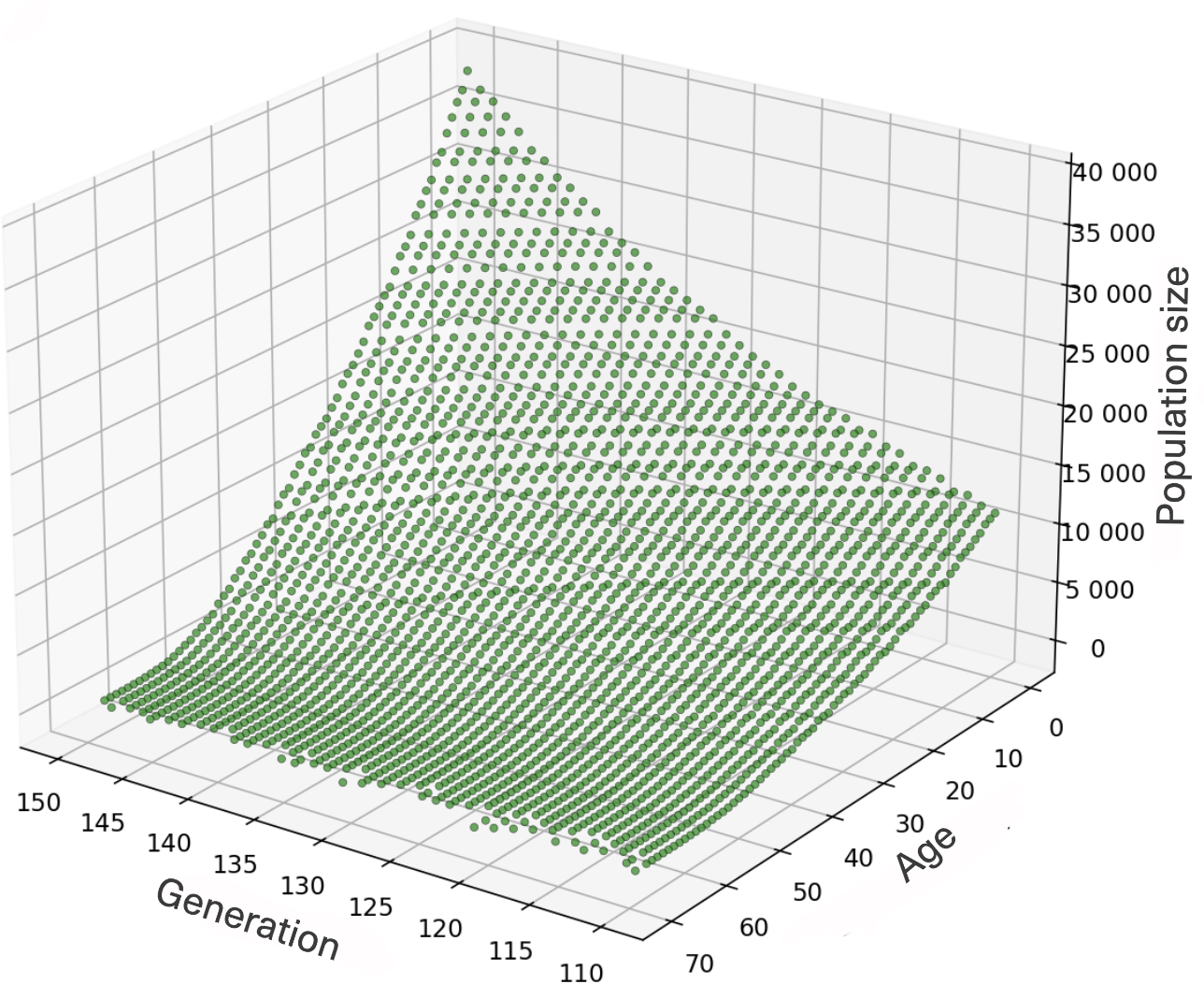}
\vspace{0.2cm}
$A = 0.001$ \\
$\alpha = 0.001$ \\
$\beta = 0.1$
\end{minipage}
\hfill
\begin{minipage}{0.3\textwidth}
\centering
\includegraphics[width=\textwidth]{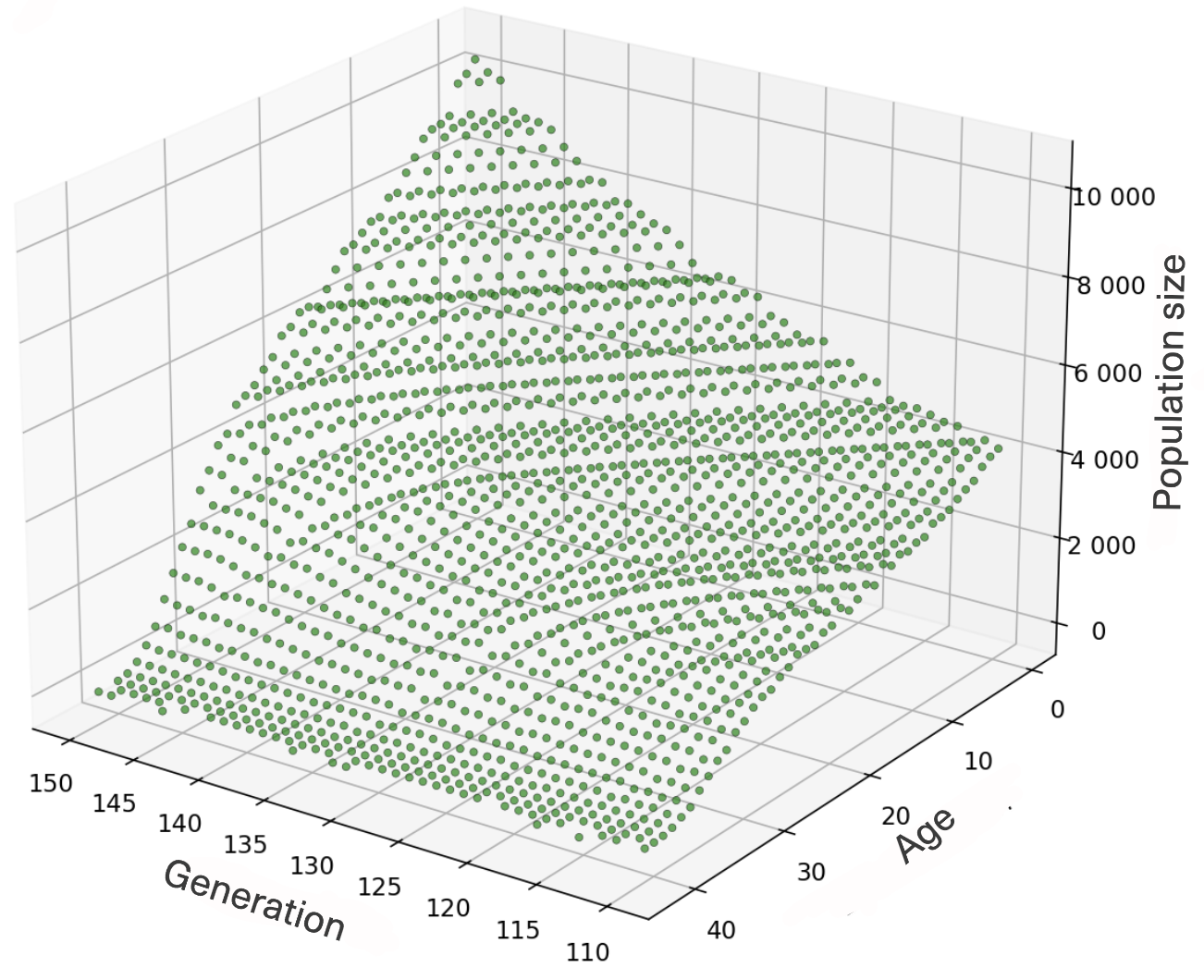}
\vspace{0.2cm}
$A = 0.001$ \\
$\alpha = 0.00001$ \\
$\beta = 0.03$
\end{minipage}
\end{center}

\noindent\textbf{Effect of fertility parameters on population dynamics}

\begin{center}
\begin{minipage}{0.3\textwidth}
\centering
\includegraphics[width=\textwidth]{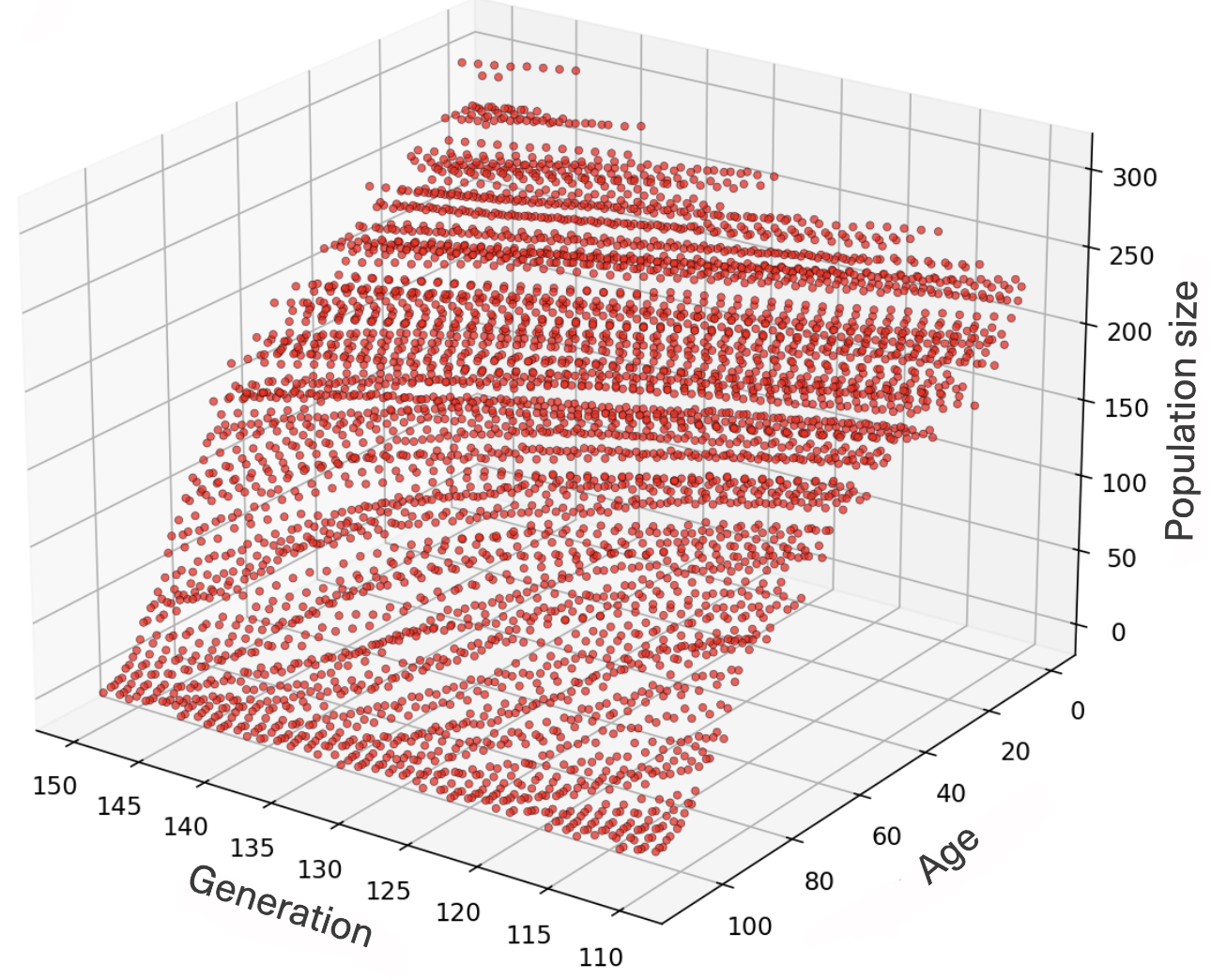}
\vspace{0.2cm}
Female individuals
\end{minipage}
\hfill
\begin{minipage}{0.3\textwidth}
\centering
\includegraphics[width=\textwidth]{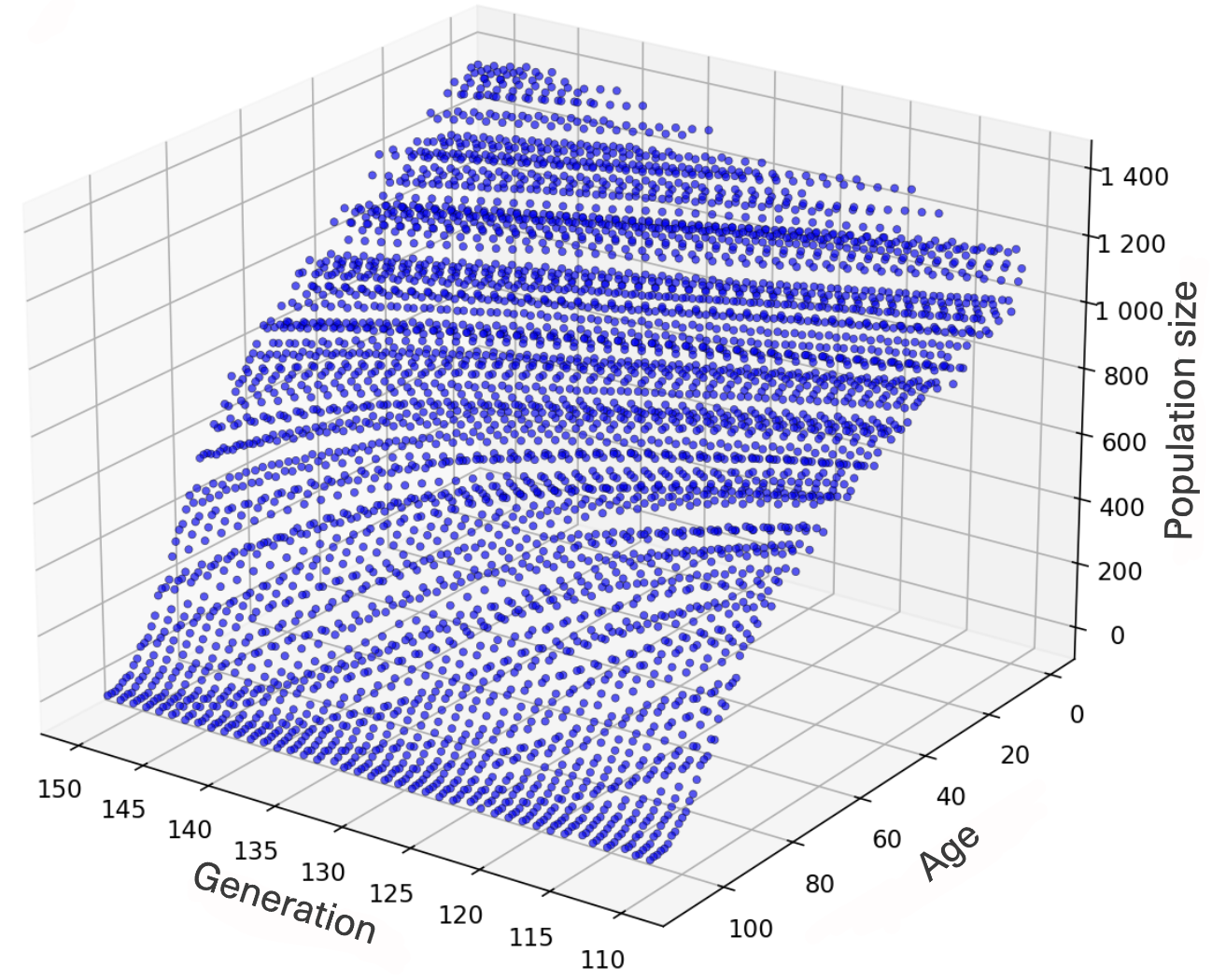}
\vspace{0.2cm}
Male individuals
\end{minipage}
\hfill
\begin{minipage}{0.3\textwidth}
\centering
\includegraphics[width=\textwidth]{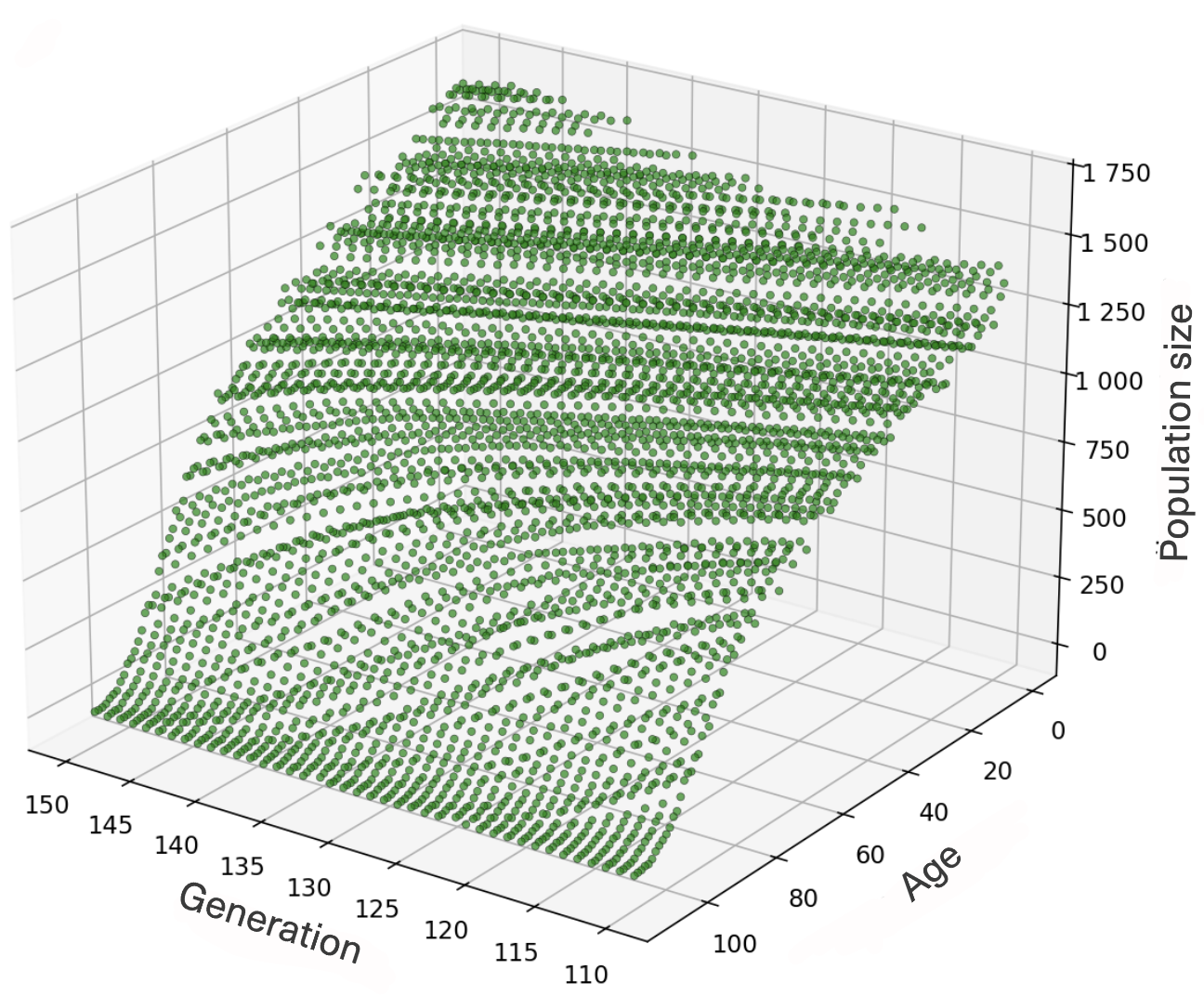}
\vspace{0.2cm}
All individuals
\end{minipage}
\vspace{0.3cm}
$p_f = 0.1\quad p_m = 0.5\quad M = 3$
\end{center}

\newpage

\section{Conclusion}

In this work, we have proposed a modification of the two-sex branching
process with random mating and overlapping generations, which differs
from the model of Molina and Mota \cite{Molina2024} by introducing
probabilistic age-dependent survival that depends on the state of the
population. For this modification, the general properties of the original
process are preserved: the zero state is absorbing, nonzero states are
transient, and the dichotomy holds (Proposition 4.1 and Theorem 4.2).

In Section 5, we established two sufficient conditions describing the
subcritical and supercritical regimes of the process. Theorem 5.1 provides
a sufficient condition for extinction: if the maximum expected number of
daughters produced by one female over her entire lifetime, taking into
account survival to each age, is less than one, then the population dies
out with probability one. Theorem 5.2 provides the reverse condition: if
the minimum expected number of daughters produced by one female over her
entire lifetime is greater than one, then the population survives with
positive probability. Together, these results generalize the theorem of
Molina and Mota \cite{Molina2024} to the case of nontrivial age-dependent
survival.

Numerical simulations of the process were performed in Python, visualizing
the dynamics of the age structure of the population for various values of
mating, survival, and fertility parameters (Section 6).

\section*{Acknowledgments}

The authors express their sincere gratitude to Professor Manuel Molina, Department of Mathematics, University of Extremadura, and Institute of Advanced Scientific Computation, University of Extremadura, Badajoz, Spain, for his valuable scientific dialogue and constructive discussions that led to the results presented in this paper. His expertise in two-sex branching processes greatly contributed to the development of this research.

\end{document}